\documentclass{elsart}
\usepackage{amsmath}
\usepackage{amssymb}
\usepackage{amsfonts}
\usepackage[all]{xy}
\newcommand{\Mat}{\mathrm{Mat}}
\newcommand{\TS}{\mathrm{TS}}

\newcommand{\K}{\mathbf{K}}

\newcommand{\gl}{\mathrm{G}}
\newcommand{\sn}{\mathfrak{S}_n}

\begin{document}

\begin{frontmatter}

\title{Linear representations, symmetric products and the commuting scheme}

\author{Francesco Vaccarino\thanksref{label1}}
\ead{francesco.vaccarino@polito.it}
\thanks[label1]{The author is partially supported
by the Research Grant 199/2004 - Politecnico di Torino}
\address{DISPEA, Politecnico di Torino, C.so Duca degli Abruzzi 24, 10129
Torino, ITALY}

\begin{abstract}
We show that the ring of multisymmetric functions over a commutative
ring is isomorphic to the ring generated by the coefficients of the characteristic polynomial
of polynomials in commuting generic matrices. As a consequence we give a surjection from the ring of invariants of several matrices to the ring of multisymmetric functions generalizing a classical result of H.Weyl and F.Junker. We also find a surjection from the ring of invariants over the commuting scheme to the ring of multisymmetric functions. This surjection is an isomophism over a characteristic zero field and induces an isomorphism at the level of reduced structures over an infinite field of positive characteristic.\end{abstract}

\end{frontmatter}
\section{Introduction}
Let $\K$ be a commutative ring.  For a $\K-$algebra $B$ we denote by $\Mat(n,B)$ be the $B-$module of $n\times n$ matrices with entries in $B$. We denote by $I_n$ the $n\times n$ identity matrix.

Let $P=\K[y_1,\dots,y_m]$, $D=\K[x_{ik}]$ and $A=\K[\xi_{kij}]$ be the polynomial rings in variables
$\{y_1,\dots,y_m\}$, $\{x_{ik}\,:\, i=1,\dots,n,\,k=1,\dots,m\}$ and  $\{\xi_{kij}\,:\, i,j=1,\dots,n,\, k=1,\dots,m\}$ over the base ring $\K$.

Following C.Procesi {\cite{dp,p1}} we introduce the generic matrices.
Let $\xi_k=(\xi_{kij})$ be the $n\times n$ matrix whose $(i,j)$ entry is $\xi_{kij}$ for $i,j=1,\dots,n$ and $k=1,\dots,m$.  We call   $\xi_{1},\dots,\xi_{m}$ the generic $n\times n$ matrices. We denote by
$A_P=\K[\xi'_{111},\dots,\xi'_{m11},\xi'_{112}\dots,\xi'_{mnn}]$
the residue algebra of $A$ modulo the ideal generated by the relations obtained from the
equation $\xi_{k}\xi_{h}=\xi_{h}\xi_{k}$ for $k,h=1,\dots,m$. Here
$\xi'_{kij}$ is the class of $\xi_{kij}$ in $A_P$. We let
$\xi'_{k}=(\xi'_{kij})$ be the $n\times n$ matrix with entries $\xi_{kij}'$ for $i,j=1,\dots,n$ and $h=1,\dots,m$. We call $\xi'_{1},\dots,\xi'_{m}$ the generic commuting $n\times n$ matrices.

There is an $n-$dimensional linear representation $\pi_P:P\to \Mat(n,A_P)$ given by mapping $y_k$ to $\xi_k'$ for $k=1,\dots,m$ (see \cite{dp}, \S1) The composition $\det\cdot\pi_P$ gives a multiplicative polynomial mapping $P\to A_P$ homogeneous of degree $n$, (see N.Bourbaki \cite{bo} A.IV.54 ).

We denote by $P^{\otimes n}$ the
tensor product $n$ times of $P$ with itself. The symmetric group
$\mathfrak{S}_n$ acts on $P^{\otimes n}$ as a group of $\K-$algebra
automorphisms by permuting the factors. We denote by $\TS^nP$ the invariants
of $P^{\otimes n}$ under $\mathfrak{S}_n$\,. By N.Roby \cite{r2} there is a unique $\K-$algebra homomorphism $\alpha:\TS^nP\to A_P$ such that $\alpha(f(y_1,\dots,y_m)^{\otimes n})=\det(\pi_P(f))=\det(f(\xi_1',\dots,\xi_m'))$. Write
\begin{equation}\det(tI_n-f(\xi_{1}',\dots,\xi_{m}'))=t^n+\sum_{k=1}^n(-1)^{k} \psi_k(f)t^{n-k}\label{psi}\end{equation}
to denote the characteristic polynomial of $\pi_P(f)=f(\xi_{1}',\dots,\xi_{m}')$.
Let $C_P$ be the subalgebra of $A_P$ generated by the coefficients of the characteristic
polynomial of $f(\xi_{1}',\dots,\xi_{m}')$ for $f\in P$.

We shall prove the following.
\begin{thm}\label{main1}
The map $\alpha:\TS^nP\to A_P$ gives an isomorphism
\[\TS^nP\cong C_P \]
i.e. the ring of symmetric tensors of order $n$ over a polynomial ring is isomorphic to the ring generated by the coefficients of the characteristic polynomial of  polynomials in generic commuting matrices.
\end{thm}

The symmetric group
$\mathfrak{S}_n$ acts on $D$ by permuting the variables so that for $\sigma\in\mathfrak{S}_n$ we have
$\sigma(x_{ik})=x_{\sigma(i)\,k}$ for all $i$ and $k$. We denote by $D^{\sn}$ the ring of the invariants for this action. It is called the ring of multisymmetric functions, see \cite{vac:fou} for a recent reference. There is an obvious $\sn-$equivariant isomorphism $D\to P^{\otimes n}$ given by mapping $x_{ik}$ to $1^{\otimes i-1}\otimes y_k\otimes 1^{\otimes n-i}$, therefore we have that $ D^{\sn}\cong\TS^nP$.

\begin{rem}
The proof of Theorem \ref{main1} will be based on Th.1 \cite{vac:fou} that gives a generating set for the ring of multisymmetric functions over a commutative ring $\K$. It can be also proved by using  results due to D.Ziplies \cite{zipgen} and F.Junkers \cite{j}.
\end{rem}

There is a surjective homomorphism of $\K-$algebra $\Delta:A\to D$ given by mapping $\xi_{kij}$ to $0$ if $i\neq j$ and to $x_{ik}$ otherwise, for $i,j=1,\dots,n$ and $k=1,\dots,m$. Observe that the $(i,i)$ entry of $\xi_{k}\xi_{h}-\xi_{h}\xi_{k}$ belongs to $\ker\Delta$ for $k,h=1,\dots,m$. Thus $\Delta$ factors through a surjective algebra homomorphism $\Delta':A_P\to D$ defined by mapping $\xi_{kij}'$ to $0$ if $i\neq j$ and to $x_{ik}$ otherwise.

The general linear group $\gl=\mathrm{GL}(n,\K)$ acts on $A$ via the action of simultaneous conjugation on $m-$tuples of matrices $\Mat(n,\K)^m$. Namely this action maps $\xi_{kij}$ to the $(i,j)$ entry of $g^{-1}\xi_kg$ for $k=1,\dots,m$ and $g\in\gl$.
Since $\xi_k\xi_h-\xi_h\xi_k=0$ is invariant by conjugation by elements of $\gl$ we have that the $\gl$ acts as a group of automorphism also on $A_P$. We denote as usual by $A_P^{\gl}$ the invariants for this action. We will prove the following Theorem.
\begin{thm}\label{main2}
The restriction of $\Delta$ to $A^{\gl}$ is a surjection onto $D^{\sn}(\cong\TS^nP)$. The same holds for the restriction of $\Delta'$ to $A_P^{\gl}$.

When $\K$ is a characteristic zero field the restriction of $\Delta'$ to $A_P^{\gl}$ gives an isomorphism $A_P^{\gl}\cong D^{\sn}$.

Let $N_P$ denote the nilradical of $A_P$.  When $\K$ is an infinite field of arbitrary characteristic $\Delta'$ induces an isomorphism between $(A_P/N_P)^{\gl}$ and $D^{\sn}$.  \end{thm}

\begin{rem}
The proof of Theorem \ref{main2} it is based on the observation that $C_P\subset A_P^{\gl}$ together with Theorem \ref{main1} . Over a characteristic zero field we have that $C_P=A_P^{\gl}$ following C.Procesi \cite{p1}, K.Sibirski{\u\i} \cite{sib}, D.Mumford \cite{mum} and the result follows.
\end{rem}

In the next section we prove Theorem \ref{main1}. The third section contains the proof of Theorem \ref{main2}, some corollaries and remarks.

\begin{ack}I would like to thank C.Procesi and M.Brion  for their hints and
discussions. A special thanks to the referee for his suggestions and patience.

This paper has developed from the notes of the talk I gave at the
AMS-Clay Summer Research Institute in Algebraic Geometry in
Seattle in July 2005: I would like to thank the organizers for
their invitation.\end{ack}

\newpage

\section{Determinant and symmetric tensors}\label{sec2}
We prove Theorem \ref{main1}.

Let $\K[x_1,\dots,x_n]$ be the polynomial ring in variables $x_1,\dots,x_n$ over $\K$. The symmetric group $\sn$ acts on $\K[x_1,\dots,x_n]$ by permuting the variables. The ring $\K[x_1,\dots,x_n]^{\sn}$ of invariants  for this action is called the ring of symmetric functions and it is generated by the elementary symmetric functions $e_1,\dots,e_n$ given by the generating function
\[\prod_{i=1}^n(t-x_i)=t^n+\sum_{k=1}^n(-1)^{k} e_kt^{n-k}\]
where the equality is calculated in $\K[t,x_1,\dots,x_n]$ with $t$ an extra variable (see \cite{md}).

Consider now an element $f\in P$, it gives a homomorphism $\rho_f:\K[x_1,\dots,x_n]\to D$ by mapping $x_i$ to $f(x_{i1},\dots,x_{im})$. We write $e_k(f)$ to denote $\rho_f(e_k)$ for $k=1,\dots,n$.
\begin{lem}\label{fou}
$D^{\sn}$ is generated by $e_1(f),\dots,e_n(f)$ with $f$ varying in the set of those monomials in $P$ that are not a proper power of another one.
\end{lem}
\begin{pf}
Th.1\cite{vac:fou}.\qed
\end{pf}
 \begin{rem}\label{ld}
 Let $\delta_k$ be the $n\times n$ diagonal matrices with $x_{1k},\dots,x_{nk}$ on the diagonal, for $k=1,\dots,m$. Consider $f(\delta_1,\dots,\delta_m)\in\Mat(n,D)$, it is a diagonal matrix having diagonal entries $\rho_f(x_1),\dots,\rho_f(x_n)$ . Hence its characteristic polynomial has $i-$th coefficient $(-1)^i\,e_i(f)$ as observed in \cite{vac:fou}.
\end{rem}

\begin{lem}
There is a unique algebra homomorphism
\[\alpha:TS^nP\to C_P\]
such that $\alpha(f\otimes\dots\otimes f)=\det(f(\xi'_1,\dots,\xi'_m))$ for
all $f\in P$. The homomorphism $\alpha$ is surjective.
\end{lem}
\begin{pf}
The homomorphism $\alpha$ is the unique $\K-$algebra homomorphism corresponding to $\det\cdot\pi_P$ which has been introduced in the paragraph just before Theorem \ref{main1} . It remains only to prove that the image of $\alpha$ is $C_P$.

Let $\alpha^*:D^{\sn}\to A_P$ be the composition of the isomorphism $D^{\sn}\cong TS^nP$ with $\alpha$. Let $t$ be an extra variable we have the homomorphism
 \[\beta=id_{\K[t]}\otimes \alpha^*:\K[t]\otimes D^{\sn}\to \K[t]\otimes A_P\] such that (recall Remark \ref{ld})
\begin{eqnarray*}
\beta(\det(t\otimes 1-1\otimes f(\delta_1,\dots,\delta_m)))&=&\beta(\prod_{i=1}^n(t\otimes 1-1\otimes\rho_f(x_{i})))\\
&=&t^n\otimes 1+\sum_{k=1}^n(-1)^{k} t^{n-k}\otimes \alpha^*(e_k(f))\\
&=& t^n\otimes 1+\sum_{k=1}^n(-1)^{k} t^{n-k}\otimes \psi_k(\pi_P(f))
\end{eqnarray*}
therefore $\alpha(\TS^nP)=C$ by Lemma \ref{fou}.\qed
\end{pf}

\begin{rem}
The generating set recalled in Lemma\ref{fou} was also known, in essence, to D.Ziplies \cite{zipgen} and F.Junkers \cite{j}.
\end{rem}

\begin{lem}\label{donto}
Let $\Delta':A_P\to D$ be given by mapping $\xi_{kij}'$ to $0$ if $i\neq j$ and to $x_{ik}$ otherwise. Then $\Delta'(C_P)=D^{\sn}$.
\end{lem}
\begin{pf}
The homomorphism $\Delta'_n:\Mat(n,A_P)\to\Mat(n,D)$ induced by $\Delta'$ is such that $\Delta'_n(\xi_k')=\delta_k$ for $k=1,\dots,m$. Thus for $f(\xi_1',\dots,\xi_m')\in \pi_P(P)\subset \Mat(n,A_P)$ we have that $\Delta_n'(f)=f(\delta_1,\dots,\delta_m)$. Hence
\begin{eqnarray*}
\det(tI_n-\Delta_n'(f))=\prod_{i=1}^n(t-\rho_f(x_{i}))&=&t^n+\sum_{k=1}^n(-1)^{k} e_k(f)t^{n-k}\\&=&t^n+\sum_{k=1}^n(-1)^{k}\Delta'( \psi_k(\pi_P(f)))t^{n-k}
\end{eqnarray*}
Thus $\Delta'(C_P)=D^{\sn}$ by Lemma \ref{fou} and the Lemma follows. \qed
\end{pf}

\begin{pf*}{Proof of Theorem \ref{main1}}
We have $\alpha^*\Delta'=id_{C_P}$ and $\Delta'\alpha^*=id_{D^{\sn}}$, thus the result follows thanks to the isomorphism $D^{\sn}\cong \TS^nP$.\qed
\end{pf*}

\section{Invariants}
We denote by $M=\Mat(n,\K)^m$ the $\K-$module of $m-$tuples of $n\times
n$ matrices.
The general linear group $\gl$ acts on $M$ by simultaneous conjugation, so that an
element $g\in\gl$ maps $(M_1,\dots,M_m)\in M$ to $(gM_1g^{-1},\dots,gM_mg^{-1})$.
This action induces another on $A$ given by mapping $\xi_{kij}$ to the $(i,j)$ entry of $g^{-1}\xi_k g$ for  $i,j=1,\dots,n$, $k=1,\dots,m$ and for all $g\in\gl$. We denote by $A^{\gl}$ the ring of invariants for this action.
Let $F=\K\{z_1,\dots,z_m\}$ be the free associative non commutative algebra on $z_1,\dots,z_m$ over the base ring $\K$. There is a linear representation $\pi_F:F\to\Mat(n,A)$ given by mapping $z_k$ to $\xi_k$ for $k=1,\dots,m$.
Consider $\pi_F(f)=f(\xi_1,\dots,\xi_m)\in\Mat(n,A)$ and let $t$ be an extra variable, we write
\begin{equation}
\det(tI_n-f)=t^n+\sum_{i=1}^n(-1)^{i}\theta_i(f)t^{n-i}
\end{equation}
We denote by $C$ the subring of $A$ generated by the coefficients $\theta_i(f)$ for $f\in F$. It is clear that $C\subset A^{\gl}$.
\begin{rem}
The previous paragraph is borrowed from \cite{dp}.
\end{rem}
Consider now the surjective homomorphism $\gamma:A\to A_P$. The kernel of $\gamma$ is the ideal of $A$ generated by the relations obtained from the equation $\xi_{k}\xi_{h}=\xi_{h}\xi_{k}$ for $k,h=1,\dots,m$ and these are invariant. Thus $A_P$ is a $\gl-$module and $\gamma$ is $\gl-$equivariant and we have a homomorphism $\gamma:A^{\gl}\to A_P^{\gl}$. Furthermore one can check that ${\gamma}(\theta_k(\pi_F(f)))=\psi_k(f(\xi_1',\dots,\xi_m'))$ for all $f\in F$ and $k=1,\dots,m$ so that ${\gamma}(C)=C_P\subset A_P^{\gl}$.

The symmetric group $\sn$ is embedded in $\gl$ via the permutation representation. Its image acts on $A$ by restriction of the previous action  by permuting the variables. Thus the ring of invariants for this action $A^{\sn}$ contains $A^{\gl}$. The kernel of $\Delta:A\to D$ is generated by the $x_{kij}$'s with $i\neq j$, therefore $\Delta$ is $\sn-$equivariant and we have a homomorphism ${\Delta}:A^{G}\to D^{\sn}$.

Clearly $A_P$ is an $\sn-$module since it is a $\gl-$module and the homomorphism $\gamma$ is $\sn-$equivariant because it is $\gl-$equivariant. Therefore $\Delta'$ is $\sn-$equivariant so that ${\Delta'}(A_P^{\gl})\subset D^{\sn}$.

\begin{lem}\label{mu}
Let $\K$ be a characteristic zero field. The general linear group $\gl$ is linear reductive thus $A^{\gl}\to A_P^{\gl}$ is surjective.
\end{lem}
\begin{pf}
\cite{mum}, Chapter\,1, pag.26.\qed
\end{pf}
\begin{lem}\label{nu}
Let $\K$ be a characteristic zero field. In this case $C=A^{\gl}$ hence $A^{\gl}_P=C_P$.
\end{lem}
\begin{pf}
C.Procesi and K.Sibirski{\u\i}\cite{p1,sib} proved that $C=A^{\gl}$. The result follows from Lemma \ref{mu}.
\end{pf}

\begin{pf*}{Proof of Theorem \ref{main2}}
Under the induced homomorphisms $\Delta_n:\Mat(n,A)\to\Mat(n,D)$ and $\Delta_n':\Mat(n,A_P)\to\Mat(n,D)$ we have that $\delta_k=\Delta_n(\xi_k)=\Delta_n'(\xi_k')$ for $k=1,\dots,m$.  This implies $\Delta(C)=\Delta'(C_P)=D^{\sn}$. Since $C\subset A^{\gl}$ and $C_P\subset A_P^{\gl} $ we have that $\Delta(A^{\gl})=\Delta'(A_P^{\gl})=D^{\sn}$ as claimed in the statement.

Suppose now $\K$ is a characteristic zero field.  We have  $\Delta':A_P^{\gl}\xrightarrow{\cong} D^{\sn}$ by Theorem \ref{main1}, Lemma \ref{mu} and Lemma \ref{nu}.

Recall that we denoted by $N_P$ the nilradical of $A_P$, it is a $\gl-$module. Hence the action of $\gl$ on $A_P$ induces another on $A_P/N_P$ so that the natural homomorphism $A_P\to A_P/N_P$ is $\gl-$equivariant. The homomorphism $\Delta'$ factors through a homomorphism $\Delta'':A_P/N_P\to D$ because $D$ is reduced. Clearly $\Delta''((A_P/N_P)^{\gl})=\Delta'(A_P^{\gl})=D^{\sn}$. 

We show now that $\Delta''(A_P/N_P)^{\gl}\xrightarrow{\cong}D^{\sn}$ when $\K$ is an infinite field of arbitrary characteristic. Let k be the algebraic closure of $\K$.  Recall that the rational  points of  $\mathrm{GL}(n,\K)$ are dense in the group $\mathrm{GL}(n, k)$ thus we can suppose $\K$ algebraically closed without any loss of generality (see \cite{pb}, \S 6.1).
Given a $m$-tuple $(Z_1,\dots,Z_m)$
of pairwise commuting matrices there is $g\in \gl$ such that
$gZ_1\,g^{-1},\dots,gZ_m\,g^{-1}$ are all in upper triangular
form.  Let then $(M_1,\dots,M_m)$ be a tuple of pairwise commuting matrices in upper triangular form.
Consider now a $1$-parameter subgroup $\lambda$ of $\gl$. We choose $\lambda$ such that
\[\lambda(t)=\left(
               \begin{array}{cccc}
                 t^{a_1} & 0 & \ldots & 0 \\
                 0 & t^{a_2} & \ldots & 0 \\
                 \vdots & \vdots & \ddots & \vdots \\
                 0 & 0 & \ldots & t^{a_n} \\
               \end{array}
             \right)
             \]
for some positive integers $a_1>a_2>\dots >a_n$.\\
For $i=1,\dots,m$ the map
\[\lambda_i:\mathbb{A}^1\backslash \{0\}\to \mathbb{A}^{n^2},\,\, t\mapsto
\lambda(t)\,M_i\,\lambda(t)^{-1}\]
can be extended to a regular map
\[\overline{\lambda}_i:\mathbb{A}^1\to\mathbb{A}^{n^2}\]
which sends the origin of $\mathbb{A}^1$ to the diagonal matrix
having the same diagonal elements as $M_i$.
It is clear that the latter belongs to the closure of the orbit of $M_i$ for
$i=1,\dots,n$.
Thus we have that in the closure of the orbit of
$(M_1,\dots,M_m)$ there is the (closed) orbit of the $m$-tuple of
diagonal matrices obtained as above. Let now $f\in (A_P/N_P)^{\gl} $ be an invariant regular function such that $\Delta''(f)=0$. Suppose $f$ is not identically zero on the orbits of tuples of commuting matrices, then there is an orbit of a tuple of diagonal matrices over which $f\neq 0$ by continuity i.e. $\Delta''(f)\neq 0$. \qed
\end{pf*}
\begin{rem}
A classical result due to F.Junker \cite{j} and restated by H.Weyl \cite{we:bo}  says that
$D^{\sn}$ is generated by the restriction of traces to the diagonal
matrices i.e. $A^{\gl}\to D^{\sn}$ is surjective over a characteristic zero field $\K$. Theorem \ref{main2} generalizes Junker - Weyl's result to any commutative base ring.
\end{rem}

\begin{cor}
Let $R\cong P/I$ be a $\K-$algebra which is a flat $\K-$module. There are two surjective homomorphisms $A^{\gl}_P\to TS^nR$ and $A_P^{\gl}\to TS^nR$. are
\label{flat}\end{cor}
\begin{pf}
Since $R$ is flat it is an inverse limit of free $\K-$modules, then
by Roby {\cite{r1}} IV, 5. Proposition IV. 5, and Bourbaki
{\cite{bo}} Exercise 8(a), AIV. p.89, we have a surjective homomorphism $\TS^nP\to\TS^nR$. The result then follows from Theorem \ref{main2}. \qed
\end{pf}

\subsection{Characteristic zero}
From now on $\K$ will be a characteristic zero field.

\begin{cor}\label{nil}
The ideal $N_P^{\gl}$ is the zero ideal and $A_P^{\gl}\cong(A_P/N_P)^{\gl}$.
\end{cor}
\begin{pf}
In characteristic zero we have $(A_P/N_P)^{\gl}\cong A_P^{\gl}/N_P^{\gl}$ since ${\gl}$ is linear reductive. We also have $A_P^{\gl}= C_P$ and
$C_P\cong D^{\sn}\cong (A_P/N_P)^{\gl}$ by Theorem \ref{main2}. Hence $N_P^{\gl}=\{0\}$ and
$A_P^{\gl}\cong(A_P/N_P)^{\gl}$.\qed
\end{pf}

\begin{rem}
Corollary \ref{nil} implies Theorem 3.3 {\cite{dom}} where $\K=\mathbb{C}$
and was also proved for
$\K=\mathbb{C}$ and $m=2$ by Gan and Ginzburg {\cite{gg}}.
\end{rem}

\begin{rem}
The affine scheme $\mathrm{Spec}\,A_P$ is called \textit{the commuting scheme}. It is conjectured that it is reduced i.e.  $N_P=0$. Corollary \ref{nil} gives some support to this conjecture.
\end{rem}

\begin{rem}
The ring $D^{S_n}$ is Cohen-Macauley by the Eagon-Hochster theorem\cite{eh}. It is also Gorenstein by Lemma 7.1.7 \cite{bk}.
Then also $A_P^{\gl}$ is Cohen-Macauley and Gorenstein by Theorem \ref{main2}.
\end{rem}

\begin{rem}
C.Procesi told me that he has independently proved the part of Theorem {\ref{main2}} relative to the characteristic zero case in this way: from H.Weyl {\cite{we:bo}} one knows that $A^{\gl}\to D^{\sn}$ is onto. Characteristic
zero implies $A^{\gl}\to A_P^{\gl}$ is onto and its kernel
contains the traces of commutators of generic matrices. Theorem
1 and Theorem 2 in {\cite{vac:fou}} jointly say that the kernel of
$A^{\gl}\to D^{\sn}$ is generated by the traces of commutators of generic matrices hence $A_P^{\gl}\cong D^{\sn}$.
\end{rem}

\end{document}